\input ./style/arxiv-general.cfg
\documentclass[aos,MSNbibl,nameyear,dvips]{arximspdf}
\makeatletter
   \@ifpackageloaded{graphicx}{}{\usepackage{graphicx}}
\makeatother
\doi{10.1214/15-AOS1368}% Updated by VTEXPTS2LaTeX.exe, 22.09.2015
\volume{44}
\issue{1}
\pubyear{2016}
\firstpage{318}
\lastpage{345}
\docsubty{FLA}

\makeatletter
\def\sfrac#1#2{#1/#2}

\def\afrac#1#2{#1/(#2)}
\def\vafrac#1#2{(#1)/(#2)}
\def\sklfrac#1#2{(#1/#2)}

\newcommand{\rrvert}{\vert}
\newcommand{\rrVert}{\Vert}
\newcommand{\llvert}{\vert}
\newcommand{\llVert}{\Vert}
\newtheorem{lemma}{Lemma}[section]
\newtheorem{proposition}{Proposition}[section]
\newtheorem{teo}{Theorem}[section]
\newtheorem{corollary}{Corollary}[section]
\newproclaim{remark}{Remark}[section]
\makeatother

\begin{document}
\begin{frontmatter}

%\dochead{}
\title{Rate exact Bayesian adaptation with modified
block~priors\thanksref{T1}}
\runtitle{Bayesian adaptive block prior}

\begin{aug}
% Corresponding author: Harrison Zhou - huibin.zhou@yale.edu% Updated by VTEXPTS2LaTeX.exe, 22.09.2015 14:57
%by VTEXPTS2LaTeX.exe, 22.09.2015 12:19
\author[A]{\fnms{Chao}~\snm{Gao}\ead[label=e1]{chao.gao@yale.edu}}
\and
\author[A]{\fnms{Harrison H.}~\snm{Zhou}\corref{}\ead[label=e2]{huibin.zhou@yale.edu}}
\runauthor{C. Gao and H.~H. Zhou}
\affiliation{Yale University}
%\dedicated{}
\address[A]{Department of Statistics\\
Yale University\\
New Haven, Connecticut 06511\\
USA\\
\printead{e1}\\
\phantom{E-mail: }\printead*{e2}}
\end{aug}
\thankstext{T1}{Supported in part by NSF Grant DMS-12-09191.}

% HISTORY:
%
\received{\smonth{8} \syear{2014}}% Updated by VTEXPTS2LaTeX.exe,
%22.09.2015 12:19
%
\revised{\smonth{7} \syear{2015}}% Updated by VTEXPTS2LaTeX.exe,
%22.09.2015 12:19

% ABSTRACT
%
\begin{abstract}
A novel block prior is proposed for adaptive Bayesian estimation. The
prior does not depend on the
smoothness of the function or the sample size. It puts sufficient prior
mass near the true signal and automatically concentrates on its
effective dimension. A rate-optimal posterior contraction is obtained
in a general framework, which includes density estimation, white noise
model, Gaussian sequence model, Gaussian regression and spectral
density estimation.
\end{abstract}

% KEYWORDS
% Pirmas kwd is didziosios raides
%
\end{frontmatter}

%s1 #&#
\section{Introduction}
\label{sec:intro}

Bayesian nonparametric estimation is attracting more and more attention
in a
wide range of applications. We consider a fundamental question in Bayesian
nonparametric estimation: is it possible to construct a prior such that the
posterior contracts to the truth with the exact optimal rate and at the same
time is adaptive regardless of the unknown smoothness? We provide a positive
answer to this question by designing a block prior on coefficients of
orthogonal series expansion of the function.

Specifically, we obtain adaptive Bayesian estimation under a Sobolev
ball assumption.
Assume that $f$ is a function on the unit interval $[0,1]$. Let $\{\phi
_{j}\}$ be the trigonometric orthogonal basis of $L^{2}[0,1]$, and
define $%
\theta_{j}=\int f\phi_{j}$ for each $j$. The Sobolev ball is
specified as
\[
E_{\alpha}(Q)= \Biggl\{ f\in L^{2}[0,1]:\sum
_{j=1}^{\infty}j^{2\alpha
}\theta_{j}^{2}
\leq Q^{2},\mbox{ with }\theta_{j}=\int f\phi
_{j}\mbox{ for each }j \Biggr\}.
\]
Under a general framework, we construct a prior $\Pi$, which satisfies the
Kullback--Leibler (KL) property and it automatically concentrates on the
effective dimension of the signal $f_{0}$, then as a consequence, the minimax
posterior contraction rate is obtained, that is,
%
%e1 #&#
\begin{equation}
P_{f_{0}}^{(n)}\Pi \bigl(\llVert f-f_{0}\rrVert
>Mn^{-\afrac{\alpha}{2\alpha
+1}}|X^{n}%
 \bigr)\longrightarrow0,
\label{postconv}
\end{equation}
where the loss function $\llVert  \cdot\rrVert  $ is the $l^{2}$-norm.

Adaptive Bayesian estimators over Sobolev balls or H\"{o}lder balls are
considered in the literature. There are two main approaches in these works.
The first one is to put a hyper-prior on the smoothness index $\alpha
$. As
is shown in \citeauthor{scricciolo06} (\citeyear{scricciolo06}) and
\citeauthor{ghosal08} (\citeyear{ghosal08}), minimax rate can be
achieved, but the set of $\alpha$ is restricted to be countable or even
finite. The second approach is to put a prior on $k$, where $k$ is the
number of basis functions for approximation, or the model dimension.
This is
called sieve prior in \citeauthor{shen01} (\citeyear{shen01}).
Examples of
using sieve prior include \citeauthor{kruijer08} (\citeyear
{kruijer08}) and %
\citeauthor{rivoirard12} (\citeyear{rivoirard12}). Their procedures are
adaptive over all $\alpha$, but the rates have extra logarithmic terms.
Other recent works in Bayesian adaptive estimation include \citeauthor
{vaart07} (\citeyear {vaart07,vaart09}), \citeauthor{jonge10} (\citeyear{jonge10}), %
\citeauthor{kruijer10} (\citeyear{kruijer10}), \citeauthor
{rousseau10} (%
\citeyear{rousseau10}), Shen, Tokdar and Ghosal (\citeyear{shen11}) and
\citeauthor {castillo14} (\citeyear {castillo14}), but the posterior
contraction rates in these works all miss a logarithmic factor.

The investigation of whether a logarithmic term is necessary in
the posterior contraction rate has fundamental implications. The
results can lead to answers to two important questions. First, is the
presence of a logarithmic term an intrinsic problem to Bayesian
adaptive nonparametric estimation? Second, is the presence of a
logarithmic term an artifact due to the current proof technique? The
answer to the first question should have an impact on statisticians'
views of the frequentist/Bayesian debate. The answer to the second
question will provide a better understanding on the famous ``prior mass
and testing'' framework [\citet{barron99};
\citet{ghosal00}] that is widely used to establish posterior
contraction results.

Compared to the previous results in the literature, the proposed block
prior is
adaptive over a continuum of smoothness, and its posterior contraction is
exactly rate-optimal. The framework for the applications of the
block prior is very general. It includes density estimation, white noise,
Gaussian sequence, regression and spectral density estimation.

At the point when the first draft of the paper was finished, we
received a manuscript by \citeauthor {hoffmann13} (\citeyear
{hoffmann13}) on Bayes adaptive estimation. They considered the similar
problem as ours and obtain the exact minimax rate by using a spike and
slab prior. However, their adaptation result for the $l_2$ loss only
holds for the white noise model. Since their proof technique takes
advantage of the Gaussian sequence structure, it cannot be immediately
extended to other model settings. In contrast, by designing a block
prior that especially works under the ``prior mass and testing''
framework, we are able to establish results for models including
density estimation, nonparametric regression and spectral density estimation.

The major difficulty of adaptation with the exact rate in various model
settings is the design of a prior distribution that satisfies the
conditions of the general prior mass and testing framework, which can
be applied to a wide range of models. This framework was pioneered by
\citeauthor{lecam73} (\citeyear{lecam73}) and \citeauthor
{schwartz65} (\citeyear{schwartz65}), and was later extended to the
nonparametric setting by \citeauthor{barron88} (\citeyear{barron88}),
\citeauthor{barron99} (\citeyear{barron99}) and Ghosal, Ghosh and van~der Vaart 
(\citeyear{ghosal00}). They proved as long as the prior
satisfies a Kullback--Leibler property and there exists a testing
procedure on the essential support of the prior, the posterior
distribution contracts to the truth with certain rate of convergence.
Though it is possible to analyze the posterior distribution according
to the Bayes formula directly as in \citeauthor {hoffmann13}
(\citeyear {hoffmann13}), the prior mass and testing framework imposes
the weakest assumption on the likelihood function, which makes it
flexible to various model settings. The price of such flexibility to
model settings is the rather strong requirements on the prior.
In our opinion, the design of a prior that satisfies the prior mass and
testing framework is the major difficulty of achieving rate-optimal
adaptation over various model settings.
The block prior we propose in this paper gives a solution to this
problem. We show that it possesses the strong properties required by
the prior mass and testing framework. Therefore, not only does it give
rate-optimal adaptation, the good posterior behavior also extends to
the settings beyond the white noise model.

The paper is organized as follows. In Section~\ref{sec:main}, we first
introduce a preliminary block prior $\bar{\Pi}$, which satisfies the
Kullback--Leibler property and concentrates on the effective dimension
of the
truth, and then we present the key result of this paper, adaptive rate-optimal
posterior contraction for a slightly modified prior $\Pi$ under a general
framework. As applications of the main results, we study adaptive Bayesian
estimation of various nonparametric models in
Section~\ref{sec:app}.
Section~\ref{sec:disc} discusses the posterior tail
probability bound and an extension of the theory to Besov balls. It
also includes discussion on why a logarithmic factor is usually present
in the Bayes nonparametric literature. The main
body of the proofs are presented in Section~\ref{sec:proof}.
Simulation and some auxiliary results of the proofs are given in the
supplement [\citet{supp}].

%s1.1 #&#
\subsection{Notations}

Throughout the paper, $\mathbb{P}$ and $\mathbb{E}$ are generic probability
and expectation operators, which are used whenever the distribution is clear
in the context. Small and big case letters denote constants which may vary
from line to line. We will not pay attention to the values of constants which
do not affect the result, unless otherwise specified. Notice these constants
may or may not be universal, which we shall make clear in the context. The
function $f$ and its Fourier coefficients $\theta=\{\theta_{j}\}$ are used
interchangeably. We say $f$ is distributed by $\Pi$ if the
corresponding $%
\theta\sim\Pi$. In the same way, the function space and the parameter
space of $f$ and $\theta$ will not be distinguished. The norm $\llVert  \cdot\rrVert  $
denotes both the $l^{2}$-norm of $f$ and the $l^{2}$-norm of $\theta$. For
two probabilities $P_{1}$ and $P_{2}$ with densities $p_{1}$ and
$p_{2}$, we
use the following divergences throughout the paper:
\begin{eqnarray*}
D(P_{1},P_{2}) &=&P_{1}\log\frac{p_{1}}{p_{2}},
\label{def:V}
\\
V(P_{1},P_{2}) &=&P_{1} \biggl(\log
\frac
{p_{1}}{p_{2}}-D(P_{1},P_{2}) \biggr)^{2},
\\
H(P_{1},P_{2}) &=& \biggl(\int (\sqrt{p_{1}}-
\sqrt{p_{2}} )^{2} \biggr)%
^{1/2}.
\end{eqnarray*}
We use $\theta_{j}$ and $\theta_{0j}$ to indicate the $j$th entries of
vectors $\theta=\{\theta_{j}\}$ and $\theta_{0}=\{\theta_{0j}\}$,
respectively. The bold notation $\bolds{\theta}_{k}$ represents the
vector $\{\theta_{j}\}_{j\in B_{k}}$ for the $k$th block. The rate $%
\varepsilon_{n}$ is always the minimax rate $\varepsilon_{n}^{2}=n^{-\afrac{
2\alpha}{2\alpha+1}}$.

%s2 #&#
\section{Main results}
\label{sec:main}

In this section, we first give some necessary background of Bayes
nonparametric estimation, then introduce a block prior and the result of
adaptive posterior contraction.

%s2.1 #&#
\subsection{Background} \label{sec:background}

Suppose we have data $X^{n}\sim P_{f_{0}}^{(n)}$, and the distribution $
P_{f_{0}}^{(n)}$ has density $p_{f_{0}}^{(n)}$ with\vspace*{1pt} respect to a dominating
measure. The posterior distribution for a prior $\Pi$ is defined to be
\[
\Pi\bigl(A|X^{n}\bigr)=\frac{\int_{A}\sklfrac
{p_{f}^{(n)}}{p_{f_{0}}^{(n)}}(X^{n})\,d\Pi
(f)}{\int\sklfrac{p_{f}^{(n)}}{p_{f_{0}}^{(n)}}(X^{n})\,d\Pi(f)}\qquad \mbox{where
}X^{n}\sim P_{f_{0}}^{(n)}.
\]
We need to bound the expectation of $\Pi (d(f,f_{0})>M\varepsilon
_{n}|X^{n}%
 )$ in this paper. To bound this quantity, it is sufficient to upper
bound the numerator and lower bound the denominator. Following %
\citeauthor{barron99} (\citeyear{barron99}) and \citeauthor
{ghosal00} (%
\citeyear{ghosal00}), this involves three steps:
\begin{longlist}[3.]
\item[1.] Show the prior $\Pi$ puts sufficient mass near the truth, that
is, we
need
\[
\Pi(K_{n})\geq\exp \bigl( -C_1n\varepsilon_n^{2}
\bigr),
\]
where $K_{n}= \{ D(P_{f_0}^{(n)},P_f^{(n)})\leq n\varepsilon_n^2,
V(P_{f_0}^{(n)},P_f^{(n)})\leq n\varepsilon_n^2  \} $.\vspace*{1pt}

\item[2.] Choose an appropriate set $\mathcal{F}_n$, and show the prior is
essentially supported on $\mathcal{F}_n$ in the sense that
\[
\Pi\bigl(\mathcal{F}_n^{c}\bigr)\leq\exp \bigl(
-C_2n\varepsilon^{2} \bigr).
\]
This controls the complexity of the prior.

\item[3.] Construct a testing function $\phi_{n}$ for the following testing
problem:
\[
H_{0}:f=f_{0}\quad\mbox{vs.}\quad H_{1}:f\in
\operatorname{supp}(\Pi )\cap\mathcal{F}_{n}\quad\mbox{and}\quad
d(f,f_{0})>M\varepsilon_{n}.
\]
The testing error needs to be well controlled in the sense that
\[
P_{f_{0}}^{(n)}\phi_{n}\vee\sup_{f\in H_{1}}P_{f}^{(n)}(1-
\phi _{n})\leq \exp \bigl( -C_{3}n\varepsilon^{2}
\bigr).
\]
\end{longlist}

Note that the constants $C_{1},C_{2} $ and $C_{3}$ are different in these
three steps above. Step~1 lower bounds the prior concentration near the
truth, which leads to a lower bound for the denominator $\int\frac
{p_{f}^{(n)}}{p_{f_{0}}^{(n)}}(X^{n})\,d\Pi(f)$. It is originated from %
\citeauthor{schwartz65} (\citeyear{schwartz65}). Steps~2~and~3 are
mainly for upper bounding the numerator $\int_A \frac
{p_{f}^{(n)}}{p_{f_{0}}^{(n)}}(X^{n})\,d\Pi(f)$. The testing idea in
step~3 is
initialized by \citeauthor{lecam73} (\citeyear{lecam73}) and %
\citeauthor{schwartz65} (\citeyear{schwartz65}). Step~2 goes back to %
\citeauthor{barron88} (\citeyear{barron88}), who proposes the idea to choose
an appropriate $\mathcal{F}_{n}$ to regularize the alternative
hypothesis in
the test, otherwise the testing function for step~3 may never exist
[see %
\citeauthor{lecam73} (\citeyear{lecam73}) and \citeauthor{barron89} (\citeyear{barron89})].

%s2.2 #&#
\subsection{The block prior \texorpdfstring{$\bar{\Pi}$}{$bar{Pi}$}}

Given a sequence $\theta=(\theta_1,\theta_2,\ldots)$ in the Hilbert
space~$l^2$. Define the blocks to be $B_k=\{l_{k},\ldots,l_{k+1}-1\}$, and $%
\{1,2,3,\ldots\}=\bigcup_{k=0}^{\infty} B_k$. Define the block size of the $k$th
block to be $n_k=l_{k+1}-l_{k}=|B_k|$. Remember the notation
$\bolds{%
\theta}_k$ represents the vector $\{\theta_j\}_{j\in B_k}$. The block
prior $%
\bar{\Pi}$ on the function $f$ is induced by a distribution on its Fourier
sequence $\{\theta_j\}$. For each $k$, let $g_k$ be a one-dimensional
density function on $\mathbb{R}^+$.

We describe $\bar{\Pi}$ as follows:
\begin{eqnarray*}
A_{k}&\sim& g_{k}\qquad\mbox{independently for each }k,
\\
\bolds{\theta}_{k}|A_{k}&\sim& N(0,A_{k}I_{n_{k}})
\qquad\mbox {independently for each }k,
\end{eqnarray*}
where $I_{n_{k}}$ is the $n_{k}\times n_{k}$ identity matrix. In this work,
we specify $l_{k}$ to be $l_{k}=[e^{k}]$. The sequence of densities $%
\{g_{k}\}$ is used to mix the scale parameter $A_{k}$ for each block,
and we
call them mixing densities. Our theory covers a class of mixing densities.
The mixing density class $\mathcal{G}$ contains all $\{g_{k}\}$ satisfying
the following properties:
\begin{longlist}[3.]
\item[1.] There exists $c_1>0$ such that, for any $k$ and $t\in
[e^{-k^2},e^{-k}] $,
%
%e2 #&#
\begin{equation}
g_k(t)\geq\exp \bigl(-c_1 e^k \bigr).
\label{mix1}
\end{equation}

\item[2.] There exists $c_2>0$, such that for any $k$,
%
%e3 #&#
\begin{equation}
\int_0^{\infty}tg_k(t)\,dt\leq4\exp
\bigl(-c_2k^2 \bigr). \label{mix2}
\end{equation}

\item[3.] There exists $c_3>0$, such that for any $k$,
%
%e4 #&#
\begin{equation}
\int_{e^{-k^2}}^{\infty}g_k(t)\,dt\leq\exp
\bigl(-c_3e^k \bigr). \label{mix3}
\end{equation}
\end{longlist}

For a function $f_0\in E_{\alpha}(Q)$, define the set
%
%e5 #&#
\begin{equation}
\mathcal{F}_n=\mathcal{F}_n(\beta)= \biggl\{\theta:\sum
_{j>(n\beta
^{-1})^{%
\afrac{1}{2\alpha+1}}}(\theta_j-\theta_{0j})^2
\leq\varepsilon _n^2 \biggr\}. \label{def:sieve}
\end{equation}
We have the following theorem characterizing the property of $\bar{\Pi}$.

%th2.1 #&#
\begin{teo}
\label{teo:main1} For the block prior $\bar{\Pi}$ with mixing
densities $%
\{g_k\}\in\mathcal{G}$, let $f_0\in E_{\alpha}(Q)$ for some $\alpha, Q>0$,
then there exists a constant $C>0$ such that
%
%e6 #&#
\begin{equation}
\bar{\Pi} \Biggl\{\sum_{j=1}^{\infty}(
\theta_j-\theta_{0j})^2\leq
\varepsilon_n^2%
 \Biggr\}\geq\exp \bigl(-Cn
\varepsilon_n^2 \bigr), \label{eq:KL}
\end{equation}
and
%
%e7 #&#
\begin{equation}
\bar{\Pi} \bigl(\mathcal{F}_n^c \bigr)\leq2\exp
\bigl(-(C+4)n\varepsilon _n^2 \bigr), \label{eq:sieve}
\end{equation}
for sufficiently large $n$ whenever $\beta\leq (\min \{
\frac{c_3}{2(C+4)},(4Q^2)^{-2\alpha} \} )%
^{2\alpha+1}$, with $c_3$ defined in (\ref{mix3}).
\end{teo}

%re2.1 #&#
\begin{remark}
The theorem presents two properties of the block prior $\bar{\Pi}$. Property
(\ref{eq:KL}) says the prior gives sufficient mass near the true
signal $%
f_{0}$. This is also recognized as the K--L condition once the
Kullback--Leibler divergence is upper bounded by the $l^{2}$-norm in the
support of the prior. Property (\ref{eq:sieve}) says the prior concentrates
on the effective dimension of the true signal $f_{0}$ automatically. In
Bayesian nonparametric theory, a testing argument is needed to prove
posterior contraction rate. Such test can be established on a sieve
receiving most of the prior mass. In (\ref{eq:sieve}), the set
$\mathcal{F}%
_{n}$ can be used as such a sieve.
\end{remark}

%re2.2 #&#
\begin{remark}
When the smoothness $\alpha$ is known, a well-known prior $\Pi
_{\alpha
}=\bigotimes_{j=1}^{\infty}N(0,j^{-2\alpha-1})$ is used in the literature.
It can be shown that this prior satisfies (\ref{eq:KL}). The block
prior $%
\bar{\Pi}$ satisfies (\ref{eq:KL}) and (\ref{eq:sieve}), and it
does not
depend on the smoothness $\alpha$. Thus, it is fully adaptive.
\end{remark}

%f1 #&#
\begin{figure}[t]

\includegraphics{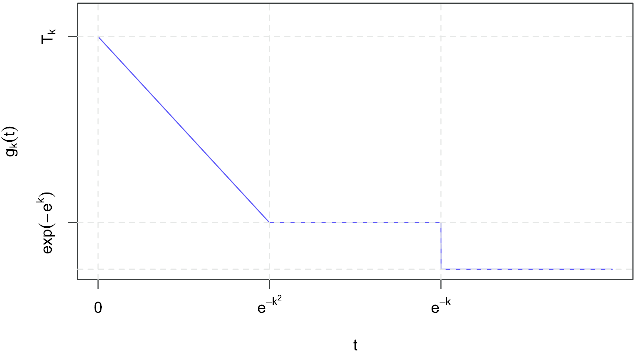}

\caption{The plot of the mixing density function $A_{k}\sim g_{k}$ defined
in (\protect\ref{prior}).}
\label{figPrior}
\end{figure}

We claim that the mixing density class $\mathcal{G}$ is not empty by
presenting an example (Figure~\ref{figPrior}):
%
%e8 #&#
\begin{equation}
g_{k}(t)=%
\cases{ e^{k^2} \bigl(\exp
\bigl(-e^k \bigr)-T_k \bigr)t+T_k, &\quad$0
\leq t\leq e^{-k^2}$;
\cr
\exp \bigl(-e^k \bigr), &
\quad$e^{-k^2}<t\leq e^{-k}$;
\cr
0, &\quad$t>
e^{-k}$.}\label{prior}
\end{equation}
The value of $T_{k}$ is specified as
%
%e9 #&#
\begin{equation}
T_{k}=2e^{k^{2}}-2\exp \bigl(-e^{k}+k^{2}-k
\bigr)+\exp \bigl(-e^{k} \bigr). \label{def:T}
\end{equation}

The following proposition is proved in the supplementary material
[\citet{supp}].

%pr2.1 #&#
\begin{proposition} \label{prop:nonempty}
The densities $\{g_k\}$ defined in (\ref{prior}) satisfies (\ref
{mix1}), (%
\ref{mix2}) and (\ref{mix3}). Thus, $\mathcal{G}$ is not empty.
\end{proposition}

%s2.3 #&#
\subsection{Adaptive posterior contraction of the modified block prior \texorpdfstring{$\Pi$}{$Pi$}}

In order to prove posterior contraction rate, it is essential to
construct a
suitable test. A~preliminary test is first constructed in a local
neighborhood. Then a global test is established by combining all the local
tests when the metric entropy is well controlled. We say the distance $d$
satisfies the testing property with respect to the prior $\Pi$ and the
truth $f_{0}$ if and only if there exists some constants $L>0$ and $\xi
\in
(0,1/2)$, such that for any $f_{1}\in\operatorname{supp}(\Pi)$
satisfying $%
d(f_{0},f_{1})>\varepsilon_{n}$, we have
%
%e10 #&#
%e11 #&#
\begin{eqnarray}
P_{f_{0}}^{(n)}\phi_{n}&\leq&\exp \bigl(-Lnd^{2}(f_{0},f_{1})
\bigr), \label{eq:test1}
\\
\sup_{\{f\in\operatorname{supp}(\Pi):d(f,f_{1})\leq\xi
d(f_{0},f_{1})\}}P_{f}^{(n)}(1-
\phi_{n})&\leq&\exp \bigl(-Lnd^{2}(f_{0},f_{1})%
 \bigr), \label{eq:test2}
\end{eqnarray}
for some testing function $\phi_{n}$. Then a global test can be
constructed for $H_{0}:f=f_{0}$ against $H_{1}= \{ f\in\mathcal
{F}%
_{n}\cap\operatorname{supp}(\Pi):d(f,f_{0})>M\varepsilon_{n} \}
$ as long as $%
d(f_{1},f_{2})\asymp\llVert  f_{1}-f_{2}\rrVert  $ for any $f_{1}$ and $f_{2}$. The
equivalence of $d$ and $\llVert  \cdot\rrVert  $ may not be true for $d$ being Hellinger
distance or total variation. We thus consider a modification of the block
prior $\bar{\Pi}$, denoted as $\Pi$, so that $d$ and $\llVert  \cdot\rrVert  $
are equivalent in
the support of the modified block prior $\Pi$. Define
\[
\Pi(A)=\frac{\bar{\Pi}(D\cap A)}{\bar{\Pi}(D)},
\]
where the constraint set $D$ needs to be designed case by case such that
\begin{eqnarray*}
D \bigl(P_{f_{1}}^{(n)},P_{f_{2}}^{(n)} \bigr)
&\leq& bn\llVert f_{1}-f_{2}\rrVert ^{2},\qquad
V \bigl(P_{f_{1}}^{(n)},P_{f_{2}}^{(n)}
\bigr)\leq bn\llVert f_{1}-f_{2}\rrVert ^{2},
\\
b^{-1}d(f_{1},f_{2})&\leq&\llVert
f_{1}-f_{2}\rrVert \leq bd(f_{1},f_{2}),
\end{eqnarray*}
for some constant $b>1$. We give a specific choice of $D$ for each model
considered in this paper. Another crucial property of $D$ we need is
that $%
\Pi$ inherits properties (\ref{eq:KL}) and (\ref{eq:sieve}) from
$\bar{\Pi}$. It is obvious that (\ref{eq:sieve}) is still true for $\Pi$ as long
as $%
\bar{\Pi}(D)>0$. Therefore, one only needs to check (\ref{eq:KL}),
which is
usually not hard as we will see in all the examples in Section~\ref{sec:app}. A general theorem covers all examples in Section~\ref{sec:app} is stated
as follows.

%th2.2 #&#
\begin{teo}
\label{teo:main2} For the block prior $\bar{\Pi}$ with mixing
densities $%
\{g_{k}\}\in\mathcal{G}$, define $\Pi(A)=\frac{\bar{\Pi}(D\cap
A)}{\bar{\Pi%
}(D)}$ with the constraint set $D$ satisfying the properties above. Let the
distance $d$ satisfy the testing property (\ref{eq:test1}) and
(\ref%
{eq:test2}). Assume that, for any $f_{0}\in E_{\alpha}(Q)\cap D$ with $
\alpha\in(\alpha^{\ast},\infty)$ and $Q\in(0,Q^{\ast})$, the
prior $%
\Pi$ inherits properties (\ref{eq:KL}) and (\ref{eq:sieve}) from
$\bar{\Pi}$
for some $C>0$. Then, for any such $f_{0}$, there exists $M>0$, such that
\[
P_{f_{0}}^{(n)}\Pi \bigl(d(f,f_{0})>Mn^{-\afrac{\alpha}{2\alpha
+1}}|%
X^{n}
\bigr)\longrightarrow0.
\]
\end{teo}

%re2.3 #&#
\begin{remark}
We note that the range $\alpha\in(\alpha^{\ast},\infty)$ and $Q\in
(0,Q^{\ast})$ is the adaptive region for the prior $\Pi$. It is determined
by the constraint set $D$ and by whether properties (\ref{eq:KL}) and
(\ref%
{eq:sieve}) can be inherited from $\bar{\Pi}$ to $\Pi$. In some examples
such as the white noise model, the modification by $D$ is not needed, so
that we have $\Pi=\bar{\Pi}$. This will result in $\alpha^{\ast
}=0$ and $%
Q^{\ast}=\infty$, and thus the prior may adapt to all Sobolev balls. In
the regression and the density estimation models, $\alpha^{\ast}$
needs to
be larger than $1/2$, and $Q^{\ast}$ can be chosen arbitrarily large by
properly picking the corresponding $D$. For the spectral density estimation,
we need $\alpha^{\ast}>3/2$. See Section~\ref{sec:app} for details.
\end{remark}

%re2.4 #&#
\begin{remark}
Theorem~\ref{teo:main2} requires the assumption $f_0\in E_{\alpha
}(Q)\cap D$. In all the nonparametric estimation examples we consider
in Section~\ref{sec:app}, we consider very specific forms of $D$ and
we are going to show that such $D$ can be removed from the assumption
because of the relation $E_{\alpha}(Q)\subset D$ for $\alpha>\alpha
^*$. This implies $E_{\alpha}(Q)\cap D=E_{\alpha}(Q)$ and we only
need $f_0\in E_{\alpha}(Q)$ in the assumption.
\end{remark}

%s3 #&#
\section{Applications}
\label{sec:app}

Given the experiment $ ( (\mathcal{X}^{(n)},\mathcal{A}%
^{(n)},P_{f}^{(n)} ):f\in E_{\alpha}(Q) )$, and observation $%
X^{n}\sim P_{f_{0}}^{(n)}$, we estimate the function $f_{0}$ by an adaptive
Bayesian procedure. The goal is to achieve the minimax posterior contraction
rate without knowing the smoothness $\alpha$. In this section, we consider
the following examples:
\begin{longlist}[2.]
\item[1. \emph{Density estimation}.] The observations $X_1,\ldots,X_n$ are i.i.d.
distributed according to the density
\[
p_f(t)=\frac{e^{f(t)}}{\int e^{f(t)}\,dt},
\]
for some function $f$ in a Sobolev ball.

\item[2. \emph{White noise}.] The observation $Y_t^{(n)}$ is from the following
process:
\[
dY_t^{(n)}=f(t)\,dt+\frac{1}{\sqrt{n}}\,dW_t,
\]
where $W_t$ is the standard Wiener process.

\item[3. \emph{Gaussian sequence}.] We have independent observations
\[
X_{i}=\theta_{i}+n^{-1/2}Z_{i},
\qquad i\in\mathbb{N},
\]
where $\{\theta_{i}\}$ are Fourier coefficients of $f$, and $\{Z_{i}\}
$ are i.i.d. standard Gaussian variables.

\item[4. \emph{Gaussian regression}.] The design is uniform $X\sim U[0,1]$.
Given $X$, $Y|X\sim N(f(X),1)$. The observations are i.i.d. pairs $%
(X_1,Y_1),\ldots,(X_n,Y_n)$.

\item[5. \emph{Spectral density}.] The observations are stationary
Gaussian time
series  $X_1,\ldots,X_n$ with mean $0$ and auto-covariance $\eta_h(g)=\int_{-%
\pi}^{\pi}e^{ih\lambda}g(\lambda)\,d\lambda$. The spectral density
$g$ is
modeled by $g=\exp (f )$ for some symmetric $f$ in a Sobolev ball.
\end{longlist}

The above models have similar frequentist estimation procedures, which is
due to the deep fact that they are asymptotically equivalent to each other
under minor regularity assumptions. References for asymptotic equivalence
theory include \citeauthor{brown96} (\citeyear{brown96}), %
\citeauthor{nussbaum96} (\citeyear{nussbaum96}), \citeauthor
{brown02} (%
\citeyear{brown02}) and %
\citeauthor{golubev10} (\citeyear{golubev10}).

%s3.1 #&#
\subsection{Density estimation}

Let $P_{f}^{(n)}$ be the product measure
$
P_{f}^{(n)}=\bigotimes_{i=1}^{n}P_{f}$.
The data is i.i.d. $X^{n}=(X_{1},\ldots,X_{n})\sim
\bigotimes_{i=1}^{n}P_{f_{0}}$. Let $P_{f}$ be dominated by Lebesgue measure
$\mu$, and it has density function
$
p_{f}(t)=\frac{e^{f(t)}}{\int_{0}^{1}e^{f(t)}\mu(dt)}$.
Consider the Fourier expansion $f=\sum_{j}\theta_{j}\phi_{j}$, and the
density $p_{f}$ can be written in the form of infinite-dimensional
exponential family:
\[
p_{f}(t)=\exp \biggl(\sum_{j}
\theta_{j}\phi_{j}(t)-\psi(\theta ) \biggr),
\]
where
\[
\psi(\theta)=\int_{0}^{1}e^{\sum_{j}\theta_{j}\phi_{j}(t)}
\mu(dt).
\]
Notice the first Fourier base function is $\phi_{1}(t)=1$. It is easy to
see that different $\theta_{1}$'s correspond to the same $p_{f}$. For
identifiability, we set $\theta_{1}=0$, so that we have $\int f(t)\mu
(dt)=\sum_{j\geq2}\theta_{j}\int\phi_{j}(t)\,dt=0$. We use the modified
block prior $\Pi(A)=\frac{\bar{\Pi}(D\cap A)}{\bar{\Pi}(D)}$ with the
constraint set
%
%e12 #&#
\begin{equation}
D= \Biggl\{ \theta:\sum_{j=1}^{\infty}\llvert
\theta_{j}\rrvert <B \Biggr\}, \label{cons:weak}
\end{equation}
for some constant $B>0$. The next lemma shows that the modified block
prior $%
\Pi$ inherits properties (\ref{eq:KL}) and (\ref{eq:sieve}) from
$\bar{\Pi}$.

%le3.1 #&#
\begin{lemma}
\label{lem:ExIn} For $\alpha^{\ast}>1/2$, define the constant
%
%e13 #&#
\begin{equation}
\gamma= \Biggl(\sum_{j=1}^{\infty}j^{-2\alpha^{\ast}}
\Biggr)^{1/2}<\infty. \label{def:const}
\end{equation}
For any $f_{0}\in E_{\alpha}(Q)$, with $\alpha\geq\alpha^{\ast}$
and $%
3\gamma Q\leq B$, there is a constant $C>0$, such that
\[
\Pi \Biggl\{ \sum_{j=1}^{\infty}(
\theta_{0j}-\theta_{j})^{2}\leq \varepsilon
_{n}^{2} \Biggr\} \geq\exp \bigl(-Cn\varepsilon_{n}^{2}
\bigr),
\]
and
\[
\Pi \bigl(\mathcal{F}_{n}^{c} \bigr)\leq2\exp
\bigl(-(C+4)n\varepsilon _{n}^{2}%
 \bigr).
\]
\end{lemma}

For density estimation, it is natural to use Hellinger distance as the
testing distance $d$. According to the testing theory in %
\citeauthor{lecam73} (\citeyear{lecam73}) and \citeauthor{ghosal00} (\citeyear{ghosal00}), it satisfies testing property (\ref%
{eq:test1}) and (\ref{eq:test2}). The next lemma establishes
equivalence among various
distances and divergences under $D$ defined in (\ref{cons:weak}).

%le3.2 #&#
\begin{lemma}
\label{lem:divequiv} On the set $D$, there exists a constant $b>1$, such
that
\begin{eqnarray*}
D(P_{f_1},P_{f_2}) &\leq& b\llVert \theta_1-
\theta_2\rrVert ^2,\qquad V(P_{f_1},P_{f_2})
\leq b\llVert \theta_1-\theta_2\rrVert ^2,
\\
b^{-1}H(P_{f_1},P_{f_2})&\leq&\llVert
\theta_1-\theta_2\rrVert \leq b H(P_{f_1},P_{f_2}).
\end{eqnarray*}
\end{lemma}

We will prove the above two lemmas in the supplementary material
[\citet{supp}]. The main
result of posterior contraction for density estimation is stated as follows.

%th3.1 #&#
\begin{teo}
\label{teo:densitypost} Let $\alpha^{\ast}>1/2$ be fixed, and
$\gamma$ is
the associated constant defined in (\ref{def:const}). For any $\alpha,Q$ satisfying $\alpha\geq\alpha^{\ast}$ and $B\geq3\gamma
Q$, there is a constant $M>0$, such that
\[
\sup_{f_0\in E_{\alpha}(Q)}P_{f_{0}}^{n}\Pi
\bigl(H(P_{f},P_{f_{0}})>M\varepsilon_{n}|X_{1},
\ldots,X_{n} \bigr)%
\longrightarrow0.
\]
\end{teo}

%re3.1 #&#
\begin{remark}
The prior $\Pi$ depends on the value of $B$, which determines the
range of
adaptation. For any $\alpha^*>1/2$ and $Q^*>0$, we can choose $B$ satisfying
$B\geq3\gamma Q^*$ ($\gamma$ depends on $\alpha^*$), such that the
prior $%
\Pi$ is adaptive for all $E_{\alpha}(Q)$ with $\alpha\geq\alpha^*$
and $%
Q\leq Q^*$.
\end{remark}

%s3.2 #&#
\subsection{White noise}

We let $P_{f}^{(n)}$ be the distribution of the following process:
\[
dY_t^{(n)}=f(t)\,dt+\frac{1}{\sqrt{n}}\,dW_t,
\qquad t\in[0,1],
\]
where $W_t$ is the standard Wiener process and the signal has Fourier
expansion $f=\sum_j\theta_j\phi_j$. This model is the simplest and most
studied nonparametric model. It is equivalent to the Gaussian sequence
model, and we have
\[
D\bigl(P_{f_0}^{(n)},P_{f}^{(n)}\bigr)=
\tfrac{1}{2}n\llVert f-f_0\rrVert ^2,\qquad V
\bigl(P_{f_0}^{(n)},P_{f}^{(n)}\bigr)=n
\llVert f-f_0\rrVert ^2.
\]
In the white noise model, it is natural to use the $l^2$ norm as the testing
distance $d$. The following lemma is from Lemma~5 in \citeauthor
{ghosal07} (%
\citeyear{ghosal07}).

%le3.3 #&#
\begin{lemma}
\label{lemWhite} Let $\phi_n= \{2%
\int(f_1(t)-f_0(t))\,dY_t^{(t)}>\llVert  f_1\rrVert  ^2-\llVert  f_0\rrVert  ^2 \}$. Then we have
\begin{eqnarray*}
P_{f_0}^{(n)}\phi_n&\leq&1-\Phi \bigl(\sqrt{n}
\llVert f_1-f_0\rrVert /2 \bigr),
\\
\sup_{\{f:\llVert  f-f_1\rrVert  \leq\llVert  f_1-f_0\rrVert  /4\}}P_{f}^{(n)}(1-
\phi_n)&\leq& 1-\Phi \bigl(%
\sqrt{n}\llVert
f_1-f_0\rrVert /4 \bigr),
\end{eqnarray*}
where $\Phi$ is the standard Gaussian cumulative distribution function.
\end{lemma}

By the property of Gaussian tail, we have
\[
1-\Phi \bigl(\sqrt{n}L\llVert  f_1-f_0\rrVert   \bigr)\leq e^{-\sklfrac{1}{2}L^2n\llVert  f_1-f_0\rrVert  ^2},
\]
provided $\sqrt{n}L\llVert  f_1-f_0\rrVert  >1$, which is true because we only need to
test those $f_1$ with $\llVert  f_1-f_0\rrVert  >M\varepsilon_n$, and we have $\sqrt{n}
\varepsilon_n\rightarrow\infty$. Therefore, in the white noise model, the
distance satisfying (\ref{eq:test1}) and (\ref{eq:test2}) is the
$l^2$ norm.
Considering that the divergence $D(P_{f_0}^{(n)},P_f^{(n)})$ and
$V(P_{f_0}^{(n)},P_f^{(n)})$ are also $l^2$ norm, we reach
the following conclusion.

%th3.2 #&#
\begin{teo}
\label{teo:whitepost} In the white noise model, for any $\alpha>0$
and $Q>0$, there exists a constant $M>0$, such
that
\[
\sup_{f_0\in E_{\alpha}(Q)}P_{f_{0}}^{(n)}\bar{\Pi} \bigl(
\llVert f-f_{0}\rrVert >M\varepsilon_{n}|Y_{t}^{(n)}
\bigr)%
\longrightarrow0.
\]
\end{teo}

Hence, this is a case that we have adaptation for all Sobolev balls.

%s3.3 #&#
\subsection{Gaussian sequence}

The Gaussian sequence model is equivalent to the while noise model. We
present this case just for illustration of the theory. Given $%
f=\sum_{j}\theta_{j}\phi_{j}$, the model $P_{f}^{(n)}$ is in a product
form
%
%e14 #&#
\begin{equation}
P_{f}^{(n)}=\bigotimes_{i=1}^{\infty}P_{\theta
_{i}}^{(n)}=
\bigotimes_{i=1}^{\infty}N\bigl(
\theta_{i},n^{-1}\bigr). \label{def:Gseq}
\end{equation}
Thus, the observations are independent Gaussian variables in the form
\[
X_{i}=\theta_{i}+n^{-1/2}Z_{i},
\qquad i\in\mathbb{N},
\]
where $\{Z_{i}\}$ are i.i.d. standard Gaussian variables. The
divergence in
this case is easy to calculate. That is,
$D(P_{f_{0}}^{(n)},P_{f}^{(n)})=\frac{n}{2}\llVert  \theta_{0}-\theta
\rrVert  ^{2}$ and $V(P_{f_{0}}^{(n)},P_{f}^{(n)})=n\llVert  \theta_{0}-\theta\rrVert  ^{2}$,
and they are exactly the $l^{2}$ norm. Define
\[
\phi_{n}(X)= \bigl\{ \llVert X-\theta_{1}\rrVert
^{2}<\llVert X-\theta_{0}\rrVert ^{2} \bigr\} =
\bigl\{ X^{T}(\theta_{1}-\theta_{0})>\llVert
\theta_{1}\rrVert ^{2}-\llVert \theta _{0}\rrVert
^{2} \bigr\}.
\]
We observe this is exactly the same test in the white noise model, and thus
Lemma~\ref{lemWhite} applies here. Therefore,
\begin{eqnarray*}
P_{f_{0}}^{(n)}\phi_{n}&\leq& e^{-\sklfrac{1}{8}n\llVert  \theta_{0}-\theta
_{1}\rrVert  ^{2}},
\\
\sup_{\{\theta:\llVert  \theta-\theta_{1}\rrVert  \leq\llVert  \theta_{1}-\theta
_{0}\rrVert  /4\}}P_{f}^{(n)}(1-
\phi_{n})&\leq& e^{-\sklfrac{1}{32}n\llVert  \theta
_{0}-\theta_{1}\rrVert  ^{2}}.
\end{eqnarray*}
The $d$ satisfying the testing property (\ref{eq:test1}) and (\ref
{eq:test2}%
) can be chosen as the $l^{2}$ norm. We thus reach the following conclusion.

%th3.3 #&#
\begin{teo}
\label{teo:seqpost} In the Gaussian sequence model, for any $\alpha
>0$ and $Q>0$, there exists a constant $%
M>0$, such that
\[
\sup_{f_0\in E_{\alpha}(Q)}P_{f_{0}}^{(n)}\bar{\Pi} \bigl(
\llVert \theta -\theta_{0}\rrVert >M\varepsilon _{n}|X_{1},X_{2},
\ldots \bigr)\longrightarrow0.
\]
\end{teo}

We have adaptation for all Sobolev balls.

%s3.4 #&#
\subsection{Gaussian regression}

We consider uniform random design instead of fixed design, because the
random design allows simple connection between various divergences and
the $%
l^{2}$ distance. The model $P_{f}^{(n)}$ gives i.i.d. observations $%
(X_{1},Y_{1}),\ldots,(X_{n},Y_{n})$ with distribution
\[
X\sim U[0,1],\qquad Y|X\sim N \bigl(f(X),1 \bigr).
\]
The theory is easily extended to general random design with $X\sim q$ for
some density $q$ on $[0,1]$ bounded from above and below. We choose the
uniform design for simplicity of presentation. The function has Fourier
expansion $f=\sum_{j}\theta_{j}\phi_{j}$ so that we can apply the modified
block prior on $f$. Let $P_{f}$ be the distribution of a single observation,
and we need to\vspace*{1pt} calculate $D(P_{f_{0}},P_{f})$ and $V(P_{f_{0}},P_{f})$.
Let $%
\phi$ be the standard normal density, and it can be shown that
$D(P_{f_{0}},P_{f})\leq\frac{1}{2}\llVert  f-f_{0}\rrVert  ^{2}$ and
$V(P_{f_{0}},P_{f})\leq (1+\frac{1}{2} (\llVert  f\rrVert  _{\infty
}^{2}+\llVert  f_{0}\rrVert  _{\infty}^{2}%
 ) )\llVert  f-f_{0}\rrVert  ^{2}$.
As what we have done in the density estimation case, we use the modified
block prior $\Pi(A)=\frac{\bar{\Pi}(A\cap D)}{\bar{\Pi}(D)}$ with the
constraint set $D= \{ \sum_{j=1}^{\infty}|\theta_{j}|<B
\} $.
According to Lemma~\ref{lem:ExIn}, the prior $\Pi$ inherits
properties (\ref%
{eq:KL}) and (\ref{eq:sieve}) from $\bar{\Pi}$. Therefore, for $f$
and $%
f_{0}\in D$,
$
V(P_{f_{0}},P_{f})\leq (1+2B^{2} )\llVert  f-f_{0}\rrVert  ^{2}$.
Next, we deal with the testing procedure. We use the likelihood ratio
test as in
the white noise and Gaussian sequence model cases, and the error is bounded
in the following lemma.

%le3.4 #&#
\begin{lemma}
\label{lem:regtest} There exists a constant $L>0$, such that for any $%
f_0,f_1\in D$ satisfying $\sqrt{n}\llVert  f_1-f_0\rrVert  >1$, there exits a testing
function $\phi_n$ with error probability bounded as
\begin{eqnarray*}
P_{f_0}^{(n)}\phi_n&\leq& e^{-Ln\llVert  f_0-f_1\rrVert  ^2},
\\
\sup_{\{f\in\operatorname{supp}(\Pi):\llVert  f-f_0\rrVert  ^2\leq
1/32\llVert  f_1-f_0\rrVert  ^2%
\}}P_{f}^{(n)}(1-
\phi_n)&\leq& e^{-Ln\llVert  f_0-f_1\rrVert  ^2}.
\end{eqnarray*}
\end{lemma}

The lemma will be proved in later sections. It says $l^{2}$ norm satisfies
the testing property (\ref{eq:test1}) and (\ref{eq:test2}). Using
Theorem~\ref{teo:main2}, we reach the following conclusion.

%th3.4 #&#
\begin{teo}
\label{teo:regpost} Let $\alpha^*>1/2$ and $\gamma$ be the constant defined
in (\ref{def:const}). In the Gaussian regression model with uniform random
design, for any $\alpha,Q$ satisfying $\alpha\geq\alpha^*$ and
$3\gamma Q\leq B$, there
exists a constant $M>0$, such that
\[
\sup_{f_0\in E_{\alpha}(Q)}P_{f_0}^{(n)}\Pi \bigl(\llVert
f-f_0\rrVert >M\varepsilon _n|X_1,
\ldots,X_n,Y_1,\ldots,Y_n \bigr)%
\longrightarrow0.
\]
\end{teo}

%re3.2 #&#
\begin{remark}
The prior $\Pi$ depends on the value of $B$, which determines the
range of
adaptation. For any $\alpha^*>1/2$ and $Q^*>0$, we can choose $B$ satisfying
$B\geq3\gamma Q^*$ ($\gamma$ depends on $\alpha^*$), such that the
prior $%
\Pi$ is adaptive for all $E_{\alpha}(Q)$ with $\alpha\geq\alpha^*$
and $%
Q\leq Q^*$.
\end{remark}

%s3.5 #&#
\subsection{Spectral density estimation}

Suppose the probability $P_{f}^{(n)}$ generates stationary Gaussian time
series data $X_1,\ldots,X_n$ with mean $0$ and spectral density $g=e^f$,
with $%
f(t)=f(-t)$. We\vspace*{1pt} assume the spectral density to be a function on $[-\pi
,\pi]$. The auto-covariance is $\eta_h=\int_{-\pi}^{\pi}e^{iht}g(t)\,dt$.
Thus, the
observation $(X_1,\ldots,X_n)$ follows $P_f^{(n)}=N (0,\Gamma
_n(g) )$,
where the covariance matrix is
\[
\Gamma_n(g)= \pmatrix{ \eta_0 & \eta_1 &
\cdots& \eta_{n-1}
\cr
\eta_1 & \eta_0 &
\cdots& \eta_{n-2}
\cr
\vdots & \vdots & \ddots& \vdots
\cr
\eta_{n-1} & \eta_{n-2} & \cdots& \eta_0}.
\]
We model the exponent of the spectral density by
$f(t)=\sum_{j=0}^{\infty}\theta_j\cos(jt)$.
According to Parseval's identity,
we have $
2\pi\llVert  g\rrVert  ^2=\llVert  \eta\rrVert  ^2$ and $2\pi\llVert  f\rrVert  ^2=\llVert  \theta\rrVert  ^2$.
We use the modified block prior $\Pi(A)=\frac{\bar{\Pi}(D\cap
A)}{\bar{\Pi}%
(D)}$ with the constraint set
%
%e15 #&#
\begin{equation}
D= \Biggl\{ \sum_{j=0}^{\infty}j\llvert
\theta_j\rrvert <B \Biggr\}. \label{cons:strong}
\end{equation}
The constraint set (\ref{cons:strong}) is stronger than (\ref{cons:weak}).
Thus, in order that the modified prior $\bar{\Pi}$ inherits
properties (\ref%
{eq:KL}) and (\ref{eq:sieve}) from the block prior $\Pi$, we need $%
\alpha>3/2 $. The following lemma will be proved in the supplementary
material [\citet{supp}].

%le3.5 #&#
\begin{lemma}
\label{lemModified3} For an arbitrary $\alpha^*>3/2$, and the
constant $%
\gamma$ defined as
%
%e16 #&#
\begin{equation}
\gamma=\sum_{j=1}^{\infty}j^{2-2\alpha^*}.
\label{def:const2}
\end{equation}
For any $f_0\in E_{\alpha}(Q)$, with $\alpha\geq\alpha^*$ and
$3\gamma Q\leq
B$, there is a constant $C>0$, such that
\[
\Pi \Biggl\{\sum_{j=1}^{\infty}(
\theta_{0j}-\theta_j)^2\leq
\varepsilon_n^2\Biggr\}\geq\exp \bigl(-Cn
\varepsilon_n^2 \bigr),
\]
and
\[
\Pi \bigl(\mathcal{F}_n^c \bigr)\leq2\exp \bigl(-(C+4)n
\varepsilon _n^2 \bigr).
\]
\end{lemma}

The\vspace*{-2pt} following lemma, comparing the $l^2$ norm with
$D(P_{f_0}^{(n)},P_f^{(n)})$ and $V(P_{f_0}^{(n)},P_{f}^{(n)})$, will be\vspace*{1pt}
proved in the supplementary material [\citet{supp}].

%le3.6 #&#
\begin{lemma}
\label{lemSpecKL} For any $f_0,f_1\in D$, we have
\begin{eqnarray*}
D\bigl(P_{f_0}^{(n)},P_{f_1}^{(n)}\bigr)
&\leq& bn\llVert f_0-f_1\rrVert ^2,
\\
V\bigl(P_{f_0}^{(n)},P_{f_1}^{(n)}\bigr)
&\leq& bn\llVert f_0-f_1\rrVert ^2,
\end{eqnarray*}
where $b>1$ is a constant only depending on $\Pi$.
\end{lemma}

The testing distance satisfying the testing properties (\ref
{eq:test1}) and (%
\ref{eq:test2}) is the \mbox{$l^2$-}norm.

%le3.7 #&#
\begin{lemma}
\label{lem:spectest} There exists constants $L>0$ and $0<\xi<1/2$,
such that
for any $f_0,f_1\in D$ with $\llVert  f_0-f_1\rrVert  ^2\geq\varepsilon_n^2$, there
exists a
testing function $\phi_n$ such that
\begin{eqnarray*}
P_{f_0}^{(n)}\phi_n&\leq&\exp \bigl(-Ln\llVert
f_0-f_1\rrVert ^2 \bigr),
\\
\sup_{\{f\in\operatorname{supp}(\Pi):\llVert  f-f_1\rrVert  \leq\xi\llVert  f_1-f_0\rrVert  \}
}P_f^{(n)}(1-%
\phi_n)&\leq&\exp \bigl(-Ln\llVert f_0-f_1
\rrVert ^2 \bigr).
\end{eqnarray*}
\end{lemma}

The lemma will be proved in later sections. We state the main result of
posterior contraction of spectral density estimation as follows.

%th3.5 #&#
\begin{teo}
\label{teo:specpost} In the spectral density estimation problem, let $%
(X_1,\ldots,\break X_n)\sim P_{f_0}^{(n)}$. For any $%
\alpha$ and $Q$ satisfying Lemma~\ref{lemModified3}, there is a
constant $%
M>0 $, such that
\[
\sup_{f_0\in E_{\alpha}(Q)}P_{f_0}^{(n)}\Pi \bigl(\llVert
f-f_0\rrVert >M\varepsilon _n|X_1,
\ldots,X_n \bigr)\longrightarrow 0.
\]
\end{teo}

%re3.3 #&#
\begin{remark}
The prior $\Pi$ depends on the value of $B$, which determines the
range of
adaptation. For any $\alpha^*>3/2$ and $Q^*>0$, we can choose $B$ satisfying
$B\geq3\gamma Q^*$ ($\gamma$ depends on $\alpha^*$), such that the
prior $%
\Pi$ is adaptive for all $E_{\alpha}(Q)$ with $\alpha\geq\alpha^*$
and $%
Q\leq Q^*$. Notice the definition of $\gamma$ in (\ref{def:const2}) is
different from that in (\ref{def:const}).
\end{remark}

%s4 #&#
\section{Discussion}\label{sec:disc}

%s4.1 #&#
\subsection{Exponential tail of the posterior}

The conclusion of the main posterior contraction result in Theorem~\ref
{teo:main2} does not specify a decaying rate of the posterior tail. In fact,
by scrutinizing the its proof, it has the following polynomial tail:
\[
P_{f_0}^{(n)}\Pi \bigl(\llVert \theta-\theta_0
\rrVert >M\varepsilon_n|X^n \bigr)\leq \frac{%
C^{\prime}}{n\varepsilon_n^2}.
\]
However, to obtain a point estimator such as posterior mean with the same
rate of convergence as $\varepsilon_n$, faster posterior tail probability is
needed [see, e.g., \citeauthor{ghosal00} (\citeyear{ghosal00}) and %
\citeauthor{shen01} (\citeyear{shen01})]. In this section, we show
that this
polynomial tail can be improved to exponential tail in all the examples we
consider in Section~\ref{sec:app}. The critical step is the following lemma,
which improves Lemma~\ref{lem:GVden} in the proof of the general
result of
Theorem~\ref{teo:main2}.

%le4.1 #&#
\begin{lemma}
\label{lem:denominator} For all statistical models we consider in
Section~\ref{sec:app} and the corresponding modified block prior $\Pi$, let
$C$ be
the constant with which $\Pi$ satisfies (\ref{eq:KL}) and (\ref{eq:sieve}).
Define
%
%e17 #&#
\begin{equation}
\mathcal{H}_n= \biggl\{\int\frac{p_f^{(n)}}{p_{f_0}^{(n)}}\bigl(X^n
\bigr)\,d\Pi (f)\geq\exp%
 \bigl(-(C+b+1)n\varepsilon_n^2
\bigr) \biggr\}. \label{def:denom}
\end{equation}
Then we have $P_{f_0}^{(n)}(\mathcal{H}_n^c)\leq\exp (-\bar
{C}n\varepsilon^2%
 )$ for $f_0\in E_{\alpha}(Q)\cap D$ and some $\bar{C}>0$.
\end{lemma}

From Lemma~\ref{lem:denominator}, we have the following improved
result for
posterior contraction.

%th4.1 #&#
\begin{teo}
\label{teo:exponentialtail} The conclusions of Theorems~\ref
{teo:densitypost}, \ref{teo:whitepost}, \ref{teo:seqpost},
\ref{teo:regpost} and~\ref{teo:specpost} can be strengthened as
\[
P_{f_0}^{(n)}\Pi \bigl(\llVert \theta-\theta_0
\rrVert >M\varepsilon_n|X^n \bigr)\leq \exp
\bigl(%
-C^{\prime}n\varepsilon_n^2 \bigr),
\]
under their corresponding settings.
\end{teo}

As a consequence, the posterior mean serves as a rate-optimal point
estimator.

%co4.1 #&#
\begin{corollary}
\label{cor:exponentialtail} Under the setting of Theorems~\ref%
{teo:densitypost}, \ref{teo:whitepost}, \ref{teo:seqpost},
\ref{teo:regpost} and~\ref{teo:specpost}, we have
\[
P_{f_0}^{(n)}\bigl\Vert \mathbb{E}_{\bar{\Pi}}\bigl(
\theta|X^n\bigr)-\theta _0\bigr\Vert
^2%
\leq M^{\prime}\varepsilon_n^2,
\]
for some constant $M^{\prime}>0$.
\end{corollary}

The proofs of Lemma~\ref{lem:denominator}, Theorem~\ref{teo:exponentialtail}
and Corollary~\ref{cor:exponentialtail} are presented in the supplementary
material [\citet{supp}].

%s4.2 #&#
\subsection{Extension to Besov balls}

Besov balls provides a more flexible collection of functions than Sobolev
balls. They are related to wavelet bases. The block prior we propose in this
paper naturally takes advantage of the multi-resolution structure of Besov
balls. Given a sequence $\{\theta_j\}$, define $\bolds{\theta}%
_k=\{\theta_{2^k+l}\}_{l=0}^{2^k-1}$ for $k=0,1,2,\ldots.$ We can view the
signals on each resolution level $\bolds{\theta}_k$ as a natural block
with size $n_k=2^k$. The Besov ball is defined as
\[
B_{p,q}^{\alpha}(Q)= \biggl\{\theta:\sum
_k 2^{skq}\llVert \bolds {\theta}%
_k
\rrVert _p^q\leq Q^q \biggr\},
\]
where $s=\alpha+\frac{1}{2}-\frac{1}{p}$ and $\llVert  \cdot\rrVert  _p$ is the
vector $%
l^p $-norm. We consider the nonsparse case where the parameters are
restricted by
%
%e18 #&#
\begin{equation}
(\alpha,p,q,Q)\in(0,\infty)\times[2,\infty]\times[1,\infty ]\times (0,\infty).
\label{eq:non-sparse}
\end{equation}
Under such restriction, the block prior is suitable for estimating the
signal in $B_{p,q}^{\alpha}(Q)$. We describe the prior $\bar{\Pi}$ as
follows:
\begin{eqnarray*}
A_k &\sim& g_k\qquad\mbox{independently for each }k,
\\
\bolds{\theta}_k|A_k&\sim& N(0,A_kI_{n_k})
\qquad\mbox {independently for each }k,
\end{eqnarray*}
where $I_{n_k}$ is the $2^k\times2^k$ identity matrix. The mixing densities
$\{g_k\}$ are defined through (\ref{prior}) and (\ref{def:T}) with the
constant $e$ replaced by $2$. It is clear that the new mixing densities
$%
\{g_k\}$ satisfies (\ref{mix1}), (\ref{mix2}) and (\ref{mix3}) with
every $e$
replaced by $2$. Define the new sieve
\[
\mathcal{F}_n= \biggl\{\sum_{k>(2\alpha+1)^{-1}\log_2(n\beta
^{-1})}
\llVert %
\bolds{\theta}_k-\bolds{\theta}_{0k}
\rrVert ^2\leq\varepsilon _n^2 \biggr\}.
\]
We state the property of the block prior $\bar{\Pi}$ targeting at Besov
balls below.

%th4.2 #&#
\begin{teo}
\label{teo:besov} For the block prior $\bar{\Pi}$ defined above, let
$%
\theta_0\in B_{p,q}^{\alpha}(Q)$ with $(\alpha,p,q,Q)$ satisfying
(\ref%
{eq:non-sparse}), then there exists a constant $C>0$ such that
%
%e19 #&#
\begin{equation}
\bar{\Pi} \Biggl\{\sum_{j=1}^{\infty}(
\theta_j-\theta_{0j})^2\leq
\varepsilon_n^2%
 \Biggr\}\geq2^{-Cn\varepsilon_n^2},
\label{eq:KLbesov}
\end{equation}
and
%
%e20 #&#
\begin{equation}
\bar{\Pi} \bigl(\mathcal{F}_n^c \bigr)
\leq2^{1-(C+4)n\varepsilon_n^2}, \label{eq:sievebesov}
\end{equation}
for sufficiently large $n$ whenever $\beta\leq (\min \{
\frac{c_3}{2(C+4)},(4Q^2)^{-2\alpha} \} )%
^{2\alpha+1}$, with $c_3$ defined in (\ref{mix3}) where $e$ is
replaced by $%
2 $.
\end{teo}

We apply the prior to the Gaussian sequence model. For other models, some
slightly extra works are needed.

%th4.3 #&#
\begin{teo}
\label{teo:white} For the Gaussian sequence model (\ref{def:Gseq})
with any $%
\theta_0\in B_{p,q}^{\alpha}(Q)$, where $(\alpha,p,q,Q)$ satisfies
(\ref%
{eq:non-sparse}), then there exists $M>0$, such that
\[
\sup_{\theta_0\in B_{p,q}^{\alpha}(Q)}P_{\theta_0}^{(n)}\bar{\Pi } \bigl(
\llVert \theta-\theta_0\rrVert > M\varepsilon_n|%
X_1,X_2,
\ldots \bigr)\longrightarrow0.
\]
Thus, the prior is adaptive for all Besov balls satisfying (\ref%
{eq:non-sparse}).
\end{teo}

We prove the results of the extension in the supplementary material
[\citet{supp}].

%s4.3 #&#
\subsection{Difficulty of achieving the exact rate}

The literature of Bayes nonparametric adaptive estimation usually
reports an extra logarithmic term along with the minimax rate $\varepsilon
_n^2$. In this section, we provide examples of two priors and
illustrate the reasons for them to have the extra logarithmic term. In
the first example, the difficulty lies in the prior itself. In the
second example, the difficulty lies in the method of proof. The
analysis also sheds light on why the block prior is able to achieve the
exact minimax rate.

%s4.3.1 #&#
\subsubsection{Difficulty due to the prior}

One of the most elegant priors on $f$ is the rescaled Gaussian process
studied by \citeauthor {vaart07} (\citeyear{vaart07,vaart09}). Consider the centered
Gaussian process $(W_t: t\in[0,1])$ with the double exponential kernel
$\mathbb{E}W_tW_s=\exp (-(s-t)^2 )$. The rescaled Gaussian
process is defined as $W_{t/c}$ for some $c$ either fixed or sampled
from a hyper-prior. The reason for the rescaling is that the original
$W_t$ has an infinitely differentiable sample path almost surely. The
rescaling step makes it rougher so that it is appropriate for
estimating a signal in Sobolev or H\"{o}lder balls. In \citeauthor
{vaart07} (\citeyear {vaart07}), the number $c$ is fixed as $
(n/(\log n)^2 )^{-1/(2\alpha+1)}$, and in \citeauthor {vaart09}
(\citeyear {vaart09}) $c$ is sampled from a Gamma distribution.\vspace*{1pt} The
posterior convergence rates are $\varepsilon_n^2 (\log n )^{\vafrac
{4\alpha}{2\alpha+1}}$ and $\varepsilon_n^2 (\log n )^{\vafrac
{4\alpha+1}{2\alpha+1}}$, respectively.

Recently, this prior was extended by \citeauthor {castillo14}
(\citeyear {castillo14}) for estimation of a function living on a
general manifold $\mathcal{M}$. They constructed a rescaled Gaussian\vspace*{1pt}
process on $\mathcal{M}$ and obtained an improved posterior
convergence rate $\varepsilon_n^2 (\log n )^{\vafrac{2\alpha
}{2\alpha+1}}$. Moreover, they also showed that such a rate cannot
further be improved by a rescaled Gaussian process with a reasonable
distribution on the rescaling parameter $c$. To be specific, they
proved that under mild conditions, there exists a function $f_0\in
B_{2,\infty}^{\alpha}(Q)$ and a constant $C>0$, such that
\[
P_{f_0}^{(n)}\Pi \bigl(\llVert f-f_0\rrVert \leq
C\varepsilon_n^2 (\log n)^{\vafrac
{2\alpha}{2\alpha+1}} |X^{n}
\bigr)\rightarrow0,
\]
for\vspace*{1pt} a rescaled Gaussian process $\Pi$. Hence, the posterior
convergence rate cannot be faster than $\varepsilon_n^2 (\log n)^{\vafrac
{2\alpha}{2\alpha+1}}$.

To summarize, in this example, the difficulty lies in the prior. It is
shown that a certain class of prior distribution is unable to achieve
the exact minimax rate.

%s4.3.2 #&#
\subsubsection{Difficulty due to the proof}

The sieve prior is another popular prior used in Bayes nonparametric
estimation. It first samples an integer $J$, which is the model
dimension. Conditioning on $J$, $\theta_j$ is sampled from some
distribution $p$ independently for all $j\leq J$ and is set to zero for
$j>J$. \citeauthor {rivoirard12} (\citeyear {rivoirard12}) considered
both fixed $J=[n^{\afrac{1}{2\alpha+1}}]$ and $J$ sampled from a
distribution with exponential tail. In the first case, the posterior
convergence rate is $\varepsilon_n^2(\log n)^2$ and a slightly slower
rate is obtained for the second case.

We argue that the difficulty for obtaining the exact minimax rate is
not due to the sieve prior itself, but due to the technique of the
proof. Using the prior mass and testing (see Section~\ref{sec:background}) proof technique developed by \citeauthor {barron99}
(\citeyear {barron99}) and Ghosal, Ghosh and van~der Vaart  (\citeyear
{ghosal00}), it is impossible to get the exact minimax rate. Let us
consider the Gaussian sequence model. In this case, the prior mass
condition for the truth $\theta_0\in E_{\alpha}(Q)$ and the rate
$\varepsilon_n^2$ is
%
%e21 #&#
\begin{equation}
\Pi \bigl(\llVert \theta-\theta_0\rrVert ^2\leq
\varepsilon_n^2 \bigr)\geq\exp \bigl(-Cn\varepsilon^2_n
\bigr),\label{eq:priormass}
\end{equation}
for some constant $C>0$. Even in the simplest sieve prior where $J$ is
chosen to be fixed, (\ref{eq:priormass}) cannot hold. This is
established in the following lemma.

%le4.2 #&#
\begin{lemma} \label{lem:counter}
Consider a sieve prior with fixed $J$ and density $p$. Assume
$\llVert  p\rrVert  _{\infty}\leq G$ for some constant $G>0$. Then,
for any $\delta_n\rightarrow0$ satisfying $\log\delta_n^{-1}\asymp
\log n$ and any $\theta_0$, we have
\[
\Pi \bigl(\llVert \theta-\theta_0\rrVert ^2\leq
\delta_n^2 \bigr)\leq\exp (-CJ\log n ),
\]
for some constant $C>0$.
\end{lemma}

In the ideal case where $J=[n^{\afrac{1}{2\alpha+1}}]$, the best
possible $\delta_n^2$ for (\ref{eq:priormass}) to hold is $\delta
_n^2\asymp n^{-\vafrac{2\alpha}{2\alpha+1}}\log n$. The extra $\log n$
term cannot be avoided to establish the desired prior mass condition.

On the other hand, we show that the sieve prior in Lemma~\ref
{lem:counter} does achieve the exact minimax rate when $p$ is taken as $N(0,1)$.

%le4.3 #&#
\begin{lemma} \label{lem:sieveprior}
For Gaussian sequence model, consider the prior distribution $\Pi
=\bigotimes_{j=1}^{J}N(0,1)$, with $J=[n^{\afrac{1}{2\alpha+1}}]$.
Then we have for any $\theta_0\in E_{\alpha}(Q)$,
\[
P_{f_0}^{(n)}\Pi \bigl(\llVert \theta-\theta_0
\rrVert ^2\geq M\varepsilon_n^2 |X^{n}
\bigr)\leq\exp \bigl(-Cn\varepsilon_n^2 \bigr),
\]
for some constants $C,M>0$.
\end{lemma}

The proof of this results takes advantage of the conjugacy and
calculates the posterior probability directly from the posterior
distribution formula. Both the proofs of Lemmas~\ref{lem:counter}~and~\ref{lem:sieveprior} are stated in the supplementary material
[\citet{supp}].

Moreover, we also establish an adaptive version of Lemma~\ref
{lem:sieveprior}. Namely, consider the prior distribution $k\sim\pi$
and conditioning on $k$, $\sqrt{n}\theta_j\sim g$ i.i.d. for $1\leq
j\leq k$ and $\theta_j=0$ for $j>k$.

%th4.4 #&#
\begin{teo}\label{teo:sieve}
Assume $\max_j\frac{\pi(j)}{\pi(j-1)}\leq c$, $-\log\pi(n^{\afrac
{1}{2\alpha+1}})\leq Cn^{\afrac{1}{2\alpha+1}}$, $|\log g(x)-\log
g(y)|\leq C(1+|x-y|)$ and $|\log g(0)|\leq C$ for some constants $c\in
(0,1)$ and $C>0$. Then, for Gaussian sequence model with any $\theta
_0\in E_{\alpha}(Q)$, we have
%
%e22 #&#
\begin{eqnarray}
\label{eq:dimension} P_{f_0}^{(n)}\Pi \bigl(k>M n^{\afrac{1}{2\alpha
+1}}
|X^{n} \bigr) &\leq& \exp \bigl(-C'n
\varepsilon_n^2 \bigr),
\nonumber\\[-8pt]\\[-8pt]\nonumber
P_{f_0}^{(n)}\Pi \bigl(\llVert \theta-
\theta_0\rrVert ^2\geq M\varepsilon _n^2
|X^{n} \bigr) &\leq& \exp \bigl(-C'n
\varepsilon_n^2 \bigr),
\end{eqnarray}
for some constants $M,C'>0$.
\end{teo}

The assumption on the prior distribution in Theorem~\ref{teo:sieve} is
mild. For example, we may choose $\pi(j)\propto e^{-Dj}$ for some
constant $D>0$ and choose $g$ to be the double exponential density. The
resulting posterior distribution contracts to the true signal at the
minimax rate adaptively for all $\alpha>0$. The success of this prior
crucially depends on the result (\ref{eq:dimension}), which allows us
to establish an optimal testing procedure on the set $J\leq Mn^{\afrac
{1}{2\alpha+1}}$. However, the proof of (\ref{eq:dimension}) takes
advantage of the independence structure of the Gaussian sequence model
and we are not able to establish (\ref{eq:dimension}) for other
models. For the same reason, the block spike and slab prior proposed in
\citeauthor {hoffmann13} (\citeyear {hoffmann13}) works only for the
Gaussian sequence model as well. Their argument in establishing (\ref
{eq:dimension}) also uses the independence structure of Gaussian
sequence model and thus does not work in other settings.

To summarize, the sieve prior is an example showing that the current
proof technique may result in the sub-optimal posterior convergence
rate, while for Gaussian sequence model, special techniques can be used
to overcome the difficulty.

%s4.3.3 #&#
\subsubsection{The block prior overcomes both difficulties}

The above discussion leads to two fundamental questions. 1. Is there a
prior which can achieve the exact minimax posterior convergence rate
without knowing $\alpha$? 2. Can the prior mass and testing proof
technique handle a minimax optimal adaptive prior? While the importance
of the first question is evident, the second question seems not that
relevant at first thought. However, the prior mass and testing method
has a great advantage that it is not specific to the choice of the
prior or the form of the model. Though we use direct calculation to
show the optimal posterior convergence in Lemma~\ref{lem:sieveprior}
and Theorem~\ref{teo:sieve}, the same proof cannot be extended to a
setting beyond Gaussian sequence model. The independence structure of
Gaussian sequence model plays an important role in the proof. In
contrast, the prior mass and testing method is very general so that it
can be applied in various settings.

The block prior provides affirmative answers to both questions. Not
only can it achieve the exact minimax rate, its proof also relies on
the prior mass and testing method, which makes it easy to apply in many
complex settings beyond Gaussian sequence model. We provide various
examples in Section~\ref{sec:app} including regression, density
estimation and spectral density estimation to illustrate the benefit of
using the prior mass and testing method. Without the prior mass and
testing method, an adaptive prior cannot be easily extended to the case
beyond Gaussian sequence model.

In fact, inequality (\ref{eq:dimension}) can be written as
%
%e23 #&#
\begin{equation}
P_{f_0}^{(n)}\Pi\bigl(\mathcal{F}_n^c|X^n
\bigr)\leq\exp\bigl(-C'n\varepsilon_n^2\bigr),
\label{eq:rev2.2}
\end{equation}
where $\mathcal{F}_n$ can be of a more general form than that in (\ref
{eq:dimension}) as long as an optimal testing procedure can be
established in $\mathcal{F}_n$. Then both the sieve prior and the
block spike and slab prior in Hoffmann, Rousseau and Schmidt-Hieber (\citeyear
{hoffmann13}) satisfy~(\ref{eq:rev2.2}). In contrast, the block prior
proposed in this paper satisfies
%
%e24 #&#
\begin{equation}
\Pi\bigl(\mathcal{F}_n^c\bigr)\leq\exp
\bigl(-C'n\varepsilon_n^2\bigr), \label{eq:rev2.3}
\end{equation}
which is one of the three conditions required by the prior mass and
testing technique. It can be shown that generally (\ref{eq:rev2.3}) is
a stronger condition than (\ref{eq:rev2.2}) in the sense that (\ref
{eq:rev2.3}) combining the prior mass lower bound imply (\ref
{eq:rev2.2}). In this sense, the block prior in this paper is a
stronger prior than the sieve prior and the block spike and slab prior
in Hoffmann, Rousseau and Schmidt-Hieber s(\citeyear {hoffmann13}). To put it in
another way, (\ref{eq:rev2.2}) is not only a condition on the prior
distribution, it is also a condition on the likelihood, which imposes
certain model structure. On the other hand, (\ref{eq:rev2.3}) is a
condition only on the prior. This is why it works in various models
besides the Gaussian sequence model.

%s5 #&#
\section{Proofs of main results}\label{sec:proof}

%s5.1 #&#
\subsection{Proof of Theorem \texorpdfstring{\protect\ref{teo:main1}}{2.1}}

\label{sec:thm:main1}

We first outline the proof and list some preparatory lemmas, and then state
the proof in detail. We introduce the notation $\bar{\Pi}^A$ to be defined
as
%
%e25 #&#
\begin{equation}
\bar{\Pi}^A=\bigotimes_{k=1}^{\infty}N(0,A_kI_{n_k}).
\label
{prior:Gaussian}
\end{equation}
Given a scale sequence $A=\{A_k\}$, the random function $f=\sum_{j}\theta_j%
\phi_j$ is distributed by $\bar{\Pi}^A$ if for each block $B_k$, $%
\bolds{\theta}_k=\{\theta_j\}_{j\in B_k}\sim N(0,A_kI_{n_k})$.
Then $%
\bar{\Pi}^A$ is a Gaussian process for a given $A$, and the block
prior is a
mixture of Gaussian process with $A$ distributed by the mixing
densities $%
\{g_k\}\in\mathcal{G}$.

Since $\bar{\Pi}$ itself is not a Gaussian process, the result for
the $l^{2}
$ small ball probability asymptotics for Gaussian process cannot be applied
directly. Our strategy is to pick a collection $V_{\alpha}$, and by
conditioning, we have
%
%e26 #&#
\begin{equation}
\bar{\Pi} (\cdot )\geq\mathbb{P}(V_{\alpha})\mathbb {E} \bigl(\bar{
\Pi%
}^{A} (\cdot ) |A\in V_{\alpha} \bigr).
\label{conditioning}
\end{equation}
Then as long as for each $A\in V_{\alpha}$, there is constants $%
C_{1},C_{2}>0$ independent of $A$, such that
%
%e27 #&#
\begin{equation}
\bar{\Pi}^{A} \Biggl\{ \sum_{j=1}^{\infty}(
\theta_{j}-\theta _{0j})^{2}\leq
\varepsilon_{n}^{2} \Biggr\} \geq\exp \bigl(-C_{1}n
\varepsilon_{n}^{2} \bigr), \label{lower1}
\end{equation}
and
%
%e28 #&#
\begin{equation}
\mathbb{P}(V_{\alpha})\geq\exp \bigl(-C_{2}n
\varepsilon_{n}^{2} \bigr), \label{lower2}
\end{equation}
then the property (\ref{eq:KL}) is a direct consequence with
$C=C_{1}+C_{2}$. Thus, picking such $V_{\alpha}$ is important. Generally speaking, for
each $A\in V_{\alpha}$, we need $\bar{\Pi}^{A}$ to behave just like a
Gaussian prior designed for estimating $f_{0}\in E_{\alpha}(Q)$ when $%
\alpha$ is known.

The distribution $\bar{\Pi}^{A}$ may be hard to deal with. Our
strategy is
to use the following simple comparison result so that we can study a simpler
distribution instead. The lemma will be proved in the supplementary
material [\citet{supp}].

%le5.1 #&#
\begin{lemma}
\label{lem:comp} For standard i.i.d. Gaussian sequence $\{Z_j\}$ and
sequences $\{a_j\}$, $\{b_j\}$ and $\{c_j\}$, suppose there is a
constant $%
R>0$ such that
\[
R^{-1}a_j\leq b_j\leq R a_j
\qquad\mbox{for all } j,
\]
then we have
\begin{eqnarray*}
\mathbb{P} \biggl(\sum_j b_j
(Z_j-c_j)^2\leq R^{-1}
\varepsilon^2 \biggr)&\leq& \mathbb{P} \biggl(\sum
_j a_j (Z_j-c_j)^2
\leq\varepsilon^2 \biggr)
\\
&\leq& \mathbb{P}%
 \biggl(\sum
_j b_j (Z_j-c_j)^2
\leq R\varepsilon^2 \biggr).
\end{eqnarray*}
\end{lemma}

Define $J_{\alpha}$ to be the smallest integer such that $J_{\alpha
}\geq
(8Q^{2})^{\afrac{1}{2\alpha}}n^{\afrac{1}{2\alpha+1}}$. Let $K$ to be the
smallest integer such that $e^{K}>J_{\alpha}$, and define $J=[e^{K}]$.
Inspired by the comparison lemma, we define
%
%e29 #&#
\begin{equation}
V_{\alpha}=V_{\alpha,R}= \biggl\{ A:R^{-1}\leq\min
_{1\leq k\leq
K}\frac{%
A_{k}}{A_{\alpha,k}}\leq\max_{1\leq k\leq K}
\frac{A_{k}}{A_{\alpha,k}}%
\leq R \biggr\}, \label{compreg}
\end{equation}
with
\[
A_{\alpha,k}=\frac{l_{k}^{-2\alpha}-l_{k+1}^{-2\alpha}}{2\alpha
(l_{k+1}-l_{k})}\qquad\mbox{for }k=1,2,\ldots,K.
\]
Define the truncated Gaussian process,
%
%e30 #&#
\begin{equation}
\bar{\Pi}_{K}^{A_{\alpha}}=\bigotimes_{k=1}^{K}N(0,A_{\alpha,k}I_{n_{k}}).
\label{prior:oracle}
\end{equation}
A random function $f=\sum_{j}\theta_{j}\phi_{j}$ is distributed by
$\bar{%
\Pi}_{K}^{A_{\alpha}}$ if $\bolds{\theta}_{k}\sim
N(0,A_{\alpha,k}I_{n_{k}})$ for each $k=1,\ldots,K$ and $\bolds{\theta}_{k}=0$
for $%
k>K $. The comparison lemma implies that we can control $\bar{\Pi
}^{A}$ for
each $A\in V_{\alpha}$ by the truncated Gaussian process $\bar{\Pi}%
_{K}^{A_{\alpha}}$. Additionally, the small ball probability of $\bar
{\Pi}%
_{K}^{A_{\alpha}}$ can be established. The argument is separated in the
following lemmas, which will be proved in later sections.

%le5.2 #&#
\begin{lemma}
\label{lem:oracle} For any $\alpha>0$, and $f_0\in E_{\alpha}(Q)$, there
exists $C_3>0$, such that
\[
\bar{\Pi}^{A_{\alpha}}_K \Biggl\{\sum
_{j=1}^{\infty}(\theta _j-
\theta_{0j})^2%
\leq\varepsilon_n^2
\Biggr\}\geq\exp \bigl(-C_3n\varepsilon_n^2
\bigr).
\]
\end{lemma}

%le5.3 #&#
\begin{lemma}
\label{lem:mixmass} For each $k$, let $A_k\sim g_k$, with $\{g_k\}\in%
\mathcal{G}$, we have
\[
\mathbb{P}(V_{\alpha})\geq\exp \bigl(-C_2n
\varepsilon_n^2 \bigr).
\]
\end{lemma}

%le5.4 #&#
\begin{lemma}
\label{lem:tailneg} For $J$ defined above, and $f_0\in E_{\alpha
}(Q)$, we
have
\[
\bar{\Pi} \biggl\{\sum_{j>J}(\theta_j-
\theta_{0j})^2\leq\frac
{\varepsilon_n^2}{2}%
 \biggr\}\geq
\frac{1}{2},
\]
for sufficiently large $n$.
\end{lemma}

\begin{pf*}{Proof of (\ref{eq:KL}) in Theorem~\ref{teo:main1}}
We first
introduce the truncated version of $\bar{\Pi}^A$ to be
\[
\bar{\Pi}^{A}_K=\bigotimes_{k=1}^KN(0,A_kI_{n_k}).
\]
By Lemma~\ref{lem:tailneg}, we have
\begin{eqnarray*}
&& \bar{\Pi} \Biggl\{\sum_{j=1}^{\infty}(
\theta_j-\theta_{0j})^2\leq
\varepsilon_n^2%
 \Biggr\}
 \\
 &&\qquad \geq \bar{\Pi} \Biggl\{
\sum_{j=1}^{J}(\theta_j-\theta
_{0j})^2\leq%
\frac{\varepsilon_n^2}{2},\sum
_{j>J}(\theta_j-\theta_{0j})^2
\leq \frac{%
\varepsilon_n^2}{2} \Biggr\}
\\
&&\qquad = \bar{\Pi} \Biggl\{\sum_{j=1}^{J}(
\theta_j-\theta_{0j})^2\leq \frac{%
\varepsilon_n^2}{2}
\Biggr\}\bar{\Pi} \biggl\{\sum_{j>J}(\theta
_j-\theta_{0j})^2%
\leq
\frac{\varepsilon_n^2}{2} \biggr\}
\\
&&\qquad \geq \frac{1}{2}\bar{\Pi} \Biggl\{\sum_{j=1}^{J}(
\theta_j-\theta _{0j})^2\leq%
\frac{\varepsilon_n^2}{2} \Biggr\},
\end{eqnarray*}
where we have used independence between different blocks in the above
equality. In the spirit of (\ref{conditioning}), we have
%
%e31 #&#
\begin{eqnarray}
&& \bar{\Pi} \Biggl\{\sum_{j=1}^{J}(
\theta_j-\theta_{0j})^2\leq\frac
{\varepsilon_n^2%
}{2}
\Biggr\}
\nonumber\\[-8pt]\\[-8pt]\nonumber
&&\qquad \geq\mathbb{P}(V_{\alpha})\mathbb{E} \Biggl(\bar{
\Pi}%
^A_K \Biggl\{\sum
_{j=1}^{\infty}(\theta_j-
\theta_{0j})^2\leq\frac
{\varepsilon_n^2%
}{2} \Biggr\} \Big| A\in
V_{\alpha} \Biggr).
\end{eqnarray}
By Lemma~\ref{lem:comp}, for each $A\in V_{\alpha}$,
\[
\bar{\Pi}^A_K \Biggl\{\sum_{j=1}^{\infty}(
\theta_j-\theta _{0j})^2\leq\frac{%
\varepsilon_n^2}{2}
\Biggr\}\geq\bar{\Pi}^{A_{\alpha}}_K \Biggl\{\sum
_{j=1}^{%
\infty}(\theta_j-
\theta_{0j})^2\leq\frac{\varepsilon_n^2}{2R} \Biggr\}.
\]
By Lemma~\ref{lem:oracle}, we have
\[
\bar{\Pi}^{A_{\alpha}}_K \Biggl\{\sum
_{j=1}^{\infty}(\theta _j-
\theta_{0j})^2%
\leq\frac{\varepsilon_n^2}{2R} \Biggr\}\geq
\exp \bigl(-C^{\prime
}n\varepsilon_n^2
\bigr).
\]
Combining what we have derived and Lemma~\ref{lem:mixmass}, (\ref
{eq:KL}) is
proved.
\end{pf*}

\begin{pf*}{Proof of (\ref{eq:sieve}) in Theorem~\ref{teo:main1}}
We
fix the
constant $C$ in (\ref{eq:KL}), and we are going to prove (\ref{eq:sieve})
with the same $C$. Remember the sieve $\mathcal{F}_{n}$ is defined by
(\ref%
{def:sieve}). Define the set
\[
\mathcal{A}_n= \biggl\{A_k\leq e^{-k^2}
\mbox{ for all }k>\frac
{1}{2\alpha+1%
}\log\bigl(n\beta^{-1}\bigr) \biggr\}.
\]
Then
\[
\bar{\Pi}\bigl(\mathcal{F}_n^c\bigr)\leq\sup
_{A\in\mathcal{A}_n}\bar{\Pi }^A\bigl(\mathcal{F%
}_n^c
\bigr)+\mathbb{P}\bigl(\mathcal{A}_n^c\bigr).
\]
Condition (\ref{mix3}) implies
\begin{eqnarray*}
\mathbb{P}\bigl(\mathcal{A}_n^c\bigr) &\leq& \sum
_{k>(2\alpha+1)^{-1}\log
(n\beta^{-1})}%
\mathbb{P} \bigl(A_k>e^{-k^2}
\bigr)
\\
&\leq& \sum_{k>(2\alpha+1)^{-1}\log(n\beta^{-1})}\exp \bigl(-c_3e^k
\bigr)
\\
&\leq& \exp \biggl(-\frac{1}{2}c_3n^{\afrac{1}{2\alpha+1}}\beta
^{-\afrac{1}{%
2\alpha+1}} \biggr)
\\
&\leq& \exp \bigl(-(C+4)n\varepsilon_n^2 \bigr).
\end{eqnarray*}
The last inequality is because $\beta\leq (\frac
{c_3}{2(C+4)} )%
^{2\alpha+1}$. We bound $\bar{\Pi}^A(\mathcal{F}_n^c)$ for each \mbox{$A\in\mathcal{A}_n$},
%
%e32 #&#
\begin{eqnarray}
\nonumber
\bar{\Pi}^A\bigl(\mathcal{F}_n^c
\bigr) &=& \bar{\Pi}^A \biggl\{\sum_{j> (n\beta^{-1})^{%
\afrac{1}{2\alpha+1}}}
(\theta_j-\theta_{0j})^2>
\varepsilon_n^2 \biggr\}
\\
\nonumber
&\leq& \bar{\Pi}^A \biggl\{2\sum
_{j> (n\beta^{-1})^{\afrac
{1}{2\alpha+1}}}\theta_j^2+2\sum
_{j> (n\beta^{-1})^{\afrac
{1}{2\alpha+1}}}\theta_{0j}^2>
\varepsilon_n^2 \biggr\}
\\
\label{eq:rev2.1} &\leq& \bar{\Pi}^A \biggl\{\sum
_{j> (n\beta
^{-1})^{\afrac{1}{2\alpha+1}}} \theta_j^2\geq
\frac{1}{4} \varepsilon_n^2 \biggr\}
\\
\nonumber
&\leq& \bar{\Pi}^A \biggl\{\sum
_{k>(2\alpha+1)^{-1}\log
(n\beta^{-1})}\llVert %
\bolds{\theta}_k\rrVert
^2\geq\frac{1}{4}\varepsilon_n^2 \biggr\}
\\
\nonumber
&\leq& \sum_{k>(2\alpha+1)^{-1}\log(n\beta^{-1})}\bar {\Pi}^A
\bigl\{\llVert %
\bolds{\theta}_k\rrVert ^2\geq
a_k\varepsilon_n^2 \bigr\},
\end{eqnarray}
where $\sum_k a_k\leq1/4$ and we choose $a_k=ak^{-2}$. The inequality
(\ref{eq:rev2.1}) is because $\theta_0\in E_{\alpha}(Q)$ and $\beta
\leq(4Q^2)^{-\vafrac{2\alpha+1}{2\alpha}}$. Define $\chi_{d}^2$ to be
the chi-square random variable with degree of freedom $d$:
\begin{eqnarray*}
&& \sum_{k>(2\alpha+1)^{-1}\log(n\beta^{-1})}\bar{\Pi}^A \bigl\{\llVert
\bolds{\theta}_k\rrVert ^2\geq a_k
\varepsilon_n^2 \bigr\}
\\
&&\qquad = \sum_{k>(2\alpha+1)^{-1}\log(n\beta^{-1})}\mathbb{P} \bigl\{
a_k^{-1}A_k%
\chi_{n_k}^2
\geq\varepsilon_n^2 \bigr\}
\\
&&\qquad = \sum_{k>(2\alpha+1)^{-1}\log(n\beta^{-1})}\mathbb{P} \bigl\{%
\varepsilon_n^{-2}C^{\prime}e^k
a_k^{-1}A_k\chi_{n_k}^2
\geq C^{\prime}e^k \bigr\}
\\
&&\qquad \leq \sum_{k>(2\alpha+1)^{-1}\log(n\beta^{-1})} \exp \bigl(-C^{\prime}e^k%
 \bigr) \bigl(1-2\varepsilon_n^{-2}C^{\prime}e^k
a_k^{-1}A_k \bigr)^{-\sfrac
{n_k}{2}},
\end{eqnarray*}
where we can choose $C^{\prime}$ sufficiently large. On the set
$\mathcal{A}%
_k$, for $n$ sufficiently large,
\[
A_k\leq e^{-k^2}\leq\frac{1}{4C^{\prime}}a_k
e^{-k}\varepsilon _n^2\qquad\mbox{for all } k>
\frac{1}{2\alpha+1}\log\bigl(n\beta^{-1}\bigr).
\]
Therefore,
\begin{eqnarray*}
&& \sum_{k>(2\alpha+1)^{-1}\log(n\beta^{-1})} \exp \bigl(-C^{\prime
}e^k
\bigr)%
 \bigl(1-2\varepsilon_n^{-2}C^{\prime}e^k
a_k^{-1}A_k \bigr)^{-\sfrac
{n_k}{2}}
\\
&&\qquad \leq \sum_{k>(2\alpha+1)^{-1}\log(n\beta^{-1})} \exp \bigl(-C^{\prime}e^k%
 \bigr) (\sqrt{2} )^{n_k}
\\
&&\qquad \leq \sum_{k>(2\alpha+1)^{-1}\log(n\beta^{-1})}\exp \biggl(-
\biggl(C^{\prime}-%
\frac{1}{2}\log2 \biggr)e^k
\biggr)
\\
&&\qquad \leq \exp \biggl(-\frac{1}{2} \biggl(C^{\prime}-
\frac{1}{2}\log 2 \biggr)\beta^{-%
\afrac{1}{2\alpha+1}}n\varepsilon^2
\biggr)
\\
&&\qquad \leq \exp \bigl(-(C+4)n\varepsilon_n^2 \bigr),
\end{eqnarray*}
with sufficiently large $C^{\prime}$ and $n$. Hence,
\[
\sup_{A\in\mathcal{A}_n}\bar{\Pi}^A\bigl(\mathcal{F}_n^c
\bigr)\leq\exp \bigl(%
-(C+4)n\varepsilon_n^2
\bigr),
\]
and we have
\[
\Pi \bigl(\mathcal{F}_n^c \bigr)\leq2\exp \bigl(-(C+4)n
\varepsilon _n^2 \bigr).
\]
Thus, the proof is complete.
\end{pf*}

%s5.2 #&#
\subsection{Proof of Theorem \texorpdfstring{\protect\ref{teo:main2}}{2.2}}\label{sec:thm:main2}

Before stating the proof of Theorem~\ref{teo:main2}, we need to establish
a testing result. It will be proved in later sections.

%le5.5 #&#
\begin{lemma}
\label{teo:test} Let $d$ be a distance satisfying the testing property
(\ref%
{eq:test1}) and (\ref{eq:test2}). Suppose that there is $b>0$ such
that for
all $f_{1},f_{2}\in D$,
\[
b^{-1}d(f_{1},f_{2})\leq\llVert
f_{1}-f_{2}\rrVert \leq bd(f_{1},f_{2}).
\]
Then for any sufficiently large $M>0$, there exists a testing function
$\phi
_{n}$, such that
\begin{eqnarray*}
P_{f_{0}}^{(n)}\phi_{n}&\leq&2\exp \bigl(-
\tfrac{1}{2}LM^{2}n\varepsilon _{n}^{2}%
 \bigr),
\\
\sup_{\{f\in\mathcal{F}_{n}\cap\operatorname{supp}(\Pi
):d(f,f_{0})>M\varepsilon
_{n}\}}P_{f}^{(n)}(1-
\phi_{n})&\leq&\exp \bigl(-L^{2}n\varepsilon _{n}^{2}
\bigr).
\end{eqnarray*}
\end{lemma}

The following result is Lemma~10 in \citeauthor{ghosal07} (%
\citeyear{ghosal07}). It lower bounds the denominator of the posterior
distribution in probability.

%le5.6 #&#
\begin{lemma}
\label{lem:GVden} Consider $\mathcal{H}_n$ defined in (\ref
{def:denom}), as
long as
\[
\Pi \bigl\{D\bigl(P_{f_0}^{(n)},P_{f}^{(n)}
\bigr)\leq bn\varepsilon_n^2, V\bigl(P_{f_0}^{(n)},P_{f}^{(n)}
\bigr)\leq bn\varepsilon_n^2 \bigr\}\geq\exp
\bigl(%
-Cn\varepsilon_n^2 \bigr),
\]
we have $P_{f_0}^{(n)}(\mathcal{H}_n^c)\leq\frac{1}{\bar
{C}^2n\varepsilon_n^2}$
for some $\bar{C}>0$.
\end{lemma}

\begin{pf*}{Proof of Theorem~\ref{teo:main2}}
Notice the prior $\Pi$
inherits the
properties (\ref{eq:KL}) and~(\ref{eq:sieve}) from $\bar{\Pi}$.
Since both $%
D(P_{f_0}^{(n)},P_{f}^{(n)})$ and $V(P_{f_0}^{(n)},P_{f}^{(n)})$ are upper
bounded by $bn\llVert  \theta_0-\theta\rrVert  ^2$, we have
\begin{eqnarray*}
&& \Pi \bigl\{D\bigl(P_{f_0}^{(n)},P_{f}^{(n)}
\bigr)\leq bn\varepsilon_n^2, V\bigl(P_{f_0}^{(n)},P_{f}^{(n)}
\bigr)\leq bn\varepsilon_n^2 \bigr\}
\\
&&\qquad \geq \Pi \Biggl\{\sum_{j=1}^{\infty}(
\theta_j-\theta_{0j})^2\leq
\varepsilon_n^2 \Biggr\} \geq\exp \bigl(-Cn
\varepsilon_n^2 \bigr),
\end{eqnarray*}
for the constant $C$ with which $\Pi$ satisfies (\ref{eq:KL}) and
(\ref%
{eq:sieve}). By Lemma~\ref{lem:GVden}, the K--L property of prior
implies $%
P_{f_0}^{(n)}(\mathcal{H}_n^c)\leq\frac{1}{\bar{C}^2n\varepsilon
_n^2}$. Let $%
\mathcal{F}_n$ be the sieve defined in~(\ref{def:sieve}) and we have
$
\Pi (\mathcal{F}_n^c )\leq2\exp (-(C+4)n\varepsilon
_n^2 )$.
Letting $\phi_n$ be the testing function in Lemma~\ref{teo:test}, we have
$P_{f_0}^{(n)}\Pi (d(f,f_0)>M\varepsilon_n|X^n )
\leq P_{f_0}^{(n)}(\mathcal{H}_n^c)+P_{f_0}^{(n)}\phi_n +
P_{f_0}^{(n)}\Pi%
 (d(f,f_0)>M\varepsilon_n|X^n )(1-\phi_n)1_{\mathcal{H}_n}$,
where the first two terms go to $0$. The last term has bound
\begin{eqnarray*}
&& P_{f_0}^{(n)}\Pi \bigl(d(f,f_0)>M
\varepsilon_n|X^n \bigr) (1-\phi _n)1_{\mathcal{H}%
_n}
\\
&&\qquad \leq \exp \bigl((C+2)n\varepsilon_n^2
\bigr)P_{f_0}^{(n)}\int_{\{f\in
\mathcal{F}%
_n:d(f,f_0)>M\varepsilon_n \}}
\frac{p_f^{(n)}}{p_{f_0}^{(n)}}%
\bigl(X^n\bigr) (1-\phi_n)
\bigl(X^n\bigr)\,d\Pi(f)
\\
&&\quad\qquad{} + \exp \bigl((C+2)n\varepsilon_n^2 \bigr)P_{f_0}^{(n)}
\int_{\mathcal
{F}_n^c}%
\frac{p_f^{(n)}}{p_{f_0}^{(n)}}\bigl(X^n
\bigr)\,d\Pi(f)
\\
&&\qquad \leq \exp \bigl((C+2)n\varepsilon_n^2 \bigr)\int
_{\{f\in\mathcal{F}%
_n:d(f,f_0)>M\varepsilon_n \}}P_{f_0}^{(n)}\frac
{p_f^{(n)}}{p_{f_0}^{(n)}}%
\bigl(X^n\bigr) (1-\phi_n) \bigl(X^n\bigr)\,d
\Pi(f)
\\
&&\qquad\quad{} + \exp \bigl((C+2)n\varepsilon_n^2 \bigr)\int
_{\mathcal
{F}_n^c}P_{f_0}^{(n)}%
\frac{p_f^{(n)}}{p_{f_0}^{(n)}}\bigl(X^n\bigr)\,d\Pi(f)
\\
&&\qquad \leq \exp \bigl((C+2)n\varepsilon_n^2 \bigr)\sup
_{\{f\in\mathcal
{F}_n\cap\operatorname{supp}(\Pi):d(f,f_0)>M\varepsilon_n\}
}P_f^{(n)}(1-\phi_n)
\\
&&\qquad\quad{} + \exp \bigl((C+2)n\varepsilon_n^2 \bigr)\Pi \bigl(
\mathcal{F}_n^c \bigr)
\\
&&\qquad\leq \exp \bigl(-\bigl(LM^2-C-2\bigr)n\varepsilon_n^2
\bigr)+2\exp \bigl(-2n\varepsilon_n^2 \bigr).
\end{eqnarray*}
We pick $M$ satisfying $M>\sqrt{L^{-1}(C+2)}$, and then every term
goes to $%
0 $. The proof is complete.
\end{pf*}

%\begin{appendix}
%\section{}
%\end{appendix}

% zodis "Acknowledgments" paliekamas pagal autoriu
\section*{Acknowledgments}
We want to thank Andrew Barron for the helpful discussion on this topic, and wish to thank
the Associate Editor and the referees for their constructive comments and
suggestions that lead to the improvement of the paper.

\begin{supplement}[id=suppA]
%\sname{Supplement A}
\stitle{Supplement to ``Rate exact Bayesian adaptation with modified block~priors''}
\slink[doi]{10.1214/15-AOS1368SUPP} %[doi,text={...}] - jei reikia
%suskaldyti doi
\sdatatype{.pdf}
\sfilename{aos1368\_supp.pdf}
\sdescription{The supplementary material [\citet{supp}] contains the remaining proofs and numerical studies of the block prior.}
\end{supplement}

% imsref loaded by linak, 2015-09-22 12:41:11

\printaddresses
\end{document}